\documentclass[10pt]{article}
\usepackage{amssymb}
\usepackage{latexsym}
\usepackage{euscript}
\usepackage{epsfig}

\newcommand{\pp}{\mathbb{P}}

\newcommand{\rr}{\mathbb{R}}

\newcommand{\zz}{\mathbb{Z}}
\newcommand{\sss}{\mathbb{S}}
\newcommand{\rot}{\mathbb{T}}

\DeclareFontFamily{OT1}{rsfs}{}
\DeclareFontShape{OT1}{rsfs}{n}{it}{<-> rsfs10}{}
\DeclareMathAlphabet{\curly}{OT1}{rsfs}{n}{it}

\newcommand{\Br}{\mathrm{Br}}
\newcommand{\Crit}{\mathrm{Crit}}

\newcommand{\Fibre}{\mathrm{Fibre}}

\newcommand{\circinf}{\sss^1_{\infty}}
\newcommand{\Homeo}{\mathrm{Homeo}}

\newcommand{\Int}{\mathrm{Int}}

\newcounter{Universal}[section]
\renewcommand{\theUniversal}{\thesection.\arabic{Universal}}

\newenvironment{Plain}{\refstepcounter{Universal} \par \vspace{0.5cm}
\noindent {\bf (\theUniversal)}\ }{\par \vspace{0.5cm}}
\newenvironment{Italic}{\refstepcounter{Universal} \par \vspace{0.5cm}
\noindent {\bf (\theUniversal)}\ \it}{\par \vspace{0.5cm}}  
  
\newenvironment{Thm}{\begin{Italic}{\sc Theorem: }}{\end{Italic}}
\newenvironment{Prop}{\begin{Italic}{\sc Proposition:} }{\end{Italic}}
\newenvironment{Defn}{\begin{Italic}{\sc Definition: } }{\end{Italic}}  
\newenvironment{Cor}{\begin{Italic}{\sc Corollary: }}{\end{Italic}}
\newenvironment{Lem}{\begin{Italic}{\sc Lemma: }}{\end{Italic}}

\newenvironment{Pf}{\par \noindent{\sc Proof:} }{\quad $\blacksquare$ \par
\vspace{0.5cm}}
 
\newenvironment{Rmk}{\begin{Plain} \noindent{\sc Remark: }}{\end{Plain}}

\newenvironment{Question}{\begin{Plain} \noindent{\sc Question:
}}{\end{Plain}}

\newcounter{enum}

\newenvironment{Eqn}{\refstepcounter{Universal} $$} {\eqno \mathrm{
(\theUniversal)} $$} 

\title{Geometric monodromy and the hyperbolic disc}
\author{Ivan Smith}
\date{New College, Oxford}

\begin{document}
\maketitle
\thispagestyle{empty}

\begin{abstract}
\noindent Symplectic four-manifolds give rise to Lefschetz fibrations,
which are determined by  monodromy representations of free groups in
mapping class groups.  We study
the topology of Lefschetz fibrations by analysing the action of the
monodromy on the universal cover of a smooth fibre and give a new and
simple proof  that Lefschetz fibrations arising from
Donaldson's 
construction via pencils of sections never decompose as non-trivial fibre
sums; in particular not all Lefschetz fibrations are fibre sums of holomorphic
Lefschetz fibrations.  We also show that there can never be
isotopy classes of simple 
closed curve invariant under the monodromy and as a corollary we give
a symplectic analogue of Manin's
theorem, showing that Lefschetz fibrations admit at most finitely many
homotopy classes of geometric section.  
\end{abstract}


\section{Introduction}

A Lefschetz fibration of a smooth four-manifold $X$ comprises a map 
$f: X \rightarrow \sss^2$ with finitely many critical points, all in distinct
fibres, at each of which $f$ takes the form $(z_1, z_2) \mapsto z_1 z_2$ with 
respect to local complex co-ordinates.  These co-ordinates (in both $X$ and 
the base sphere) are obliged to respect global orientations.  Away from the
critical values of the map $f$ we have a locally trivial fibre bundle
with fibre a smooth surface of some genus $g$;  at the finitely many critical
values the fibre has a single positive node.  The vanishing
cycles of a Lefschetz fibration are simple closed curves in a fixed fibre 
which shrink to the nodal points of the critical fibres; they are uniquely 
defined to isotopy given a choice of paths in $\sss^2$ from the fixed 
base-point to the
points of the critical set $\{ \Crit \}$.  The monodromy homomorphism
$\rho: \pi_1 (\sss^2 \backslash \{ \Crit \} ) \rightarrow \Gamma_g$ takes the
simple closed curve about one critical point to the positive Dehn twist
about the corresponding vanishing cycle.  

\vspace{0.2cm}

\noindent It follows, by an easy
adaptation of an argument due to Thurston (e.g. \cite{GompfS}, 
\cite{ivanthesis}),  that
the four-manifold $X$ is symplectic with $f$ giving rise to a distinguished 
deformation class of symplectic form, with respect to which all the fibres are
symplectic submanifolds away from their singular points.  A resurgence of 
interest in manifolds 
with this particular structure follows a remarkable theorem of Donaldson 
\cite{skd}:  any symplectic four-manifold admits such a topological structure
after finitely many symplectic blow-ups, and with the fibres being
symplectic submanifolds.  More particularly, Donaldson proves that certain
``pencils'' of sections of a line bundle with first Chern class represented
by a large multiple of the symplectic form (perturbed and scaled to be
an integral cohomology class) give rise to families of symplectic submanifolds
all of which intersect transversely pairwise at finitely many base-points.
The picture is familiar from pencils of divisors in K\"ahler surfaces, and 
as in that situation one obtains a Lefschetz fibration on blowing up these
base-points.  This material, and more, may be found in
\cite{ivanthesis} or \cite{ivanhodge} for 
instance.

\vspace{0.2cm}

\noindent It follows that the symplectic Lefschetz fibrations provided by 
Donaldson's existence theorem have some additional structure: a number of
sections of the map $f$ each of which has image a symplectic sphere of 
square $(-1)$.  The first result of this paper indicates the constraint
that this imposes on the rest of the topology of the fibration.  Recall that
given two Lefschetz fibrations $X, Y$ with fibres of the same genus,
one may form
their \emph{fibre sum};  topologically we remove neighbourhoods of
smooth fibres 
in each and glue the resulting (trivial) surface bundles over circles by any 
fibre-preserving diffeomorphism $\phi$.  Denote the result $X
\sharp_{\phi} Y$.   Symplectically we may perform the same
operation (perhaps after a suitable scale change to make the distinguished 
smooth fibres have the same symplectic area).  The fibre summation is said
to be non-trivial provided that each of the original fibrations has a non-zero
number of critical fibres (equivalently, provided the final four-manifold does
not have Euler characteristic equal to that of either of the constituent 
pieces). 

\begin{Thm} \label{withthedisc}
If $f: X \rightarrow \sss^2$ contains a section of square $(-1)$ then it 
cannot be decomposed as any non-trivial fibre sum.
\end{Thm}

\noindent Gang Tian asked if simply connected Lefschetz fibrations
are necessarily fibre sums of holomorphic fibrations.  Rather special
examples due to
the author (presented later) showed this to be false in general.
These examples relied on deep facts on symplectic manifolds arising
from Seiberg-Witten theory.  Independently, Stipsicz
\cite{Stipsicz} proved the above theorem also using Seiberg-Witten
techniques, and in particular the fact that a symplectic four-manifold
with $b_+ > 1$ cannot contain a symplectic sphere of non-negative
square.  By contrast, the proof of (\ref{withthedisc}) that
we shall present in the third section is both
self-contained and elementary.  Note in any case that Donaldson's
existence result for pencils and (\ref{withthedisc})
imply a negative answer to Tian's question:

\begin{Cor}
There are infinitely many simply connected symplectic
Lefschetz fibrations which are not fibre sums of holomorphic fibrations.
\end{Cor}

\noindent Our method provides a nice characterisation of the self-intersection
of a section of a Lefschetz fibration in terms of a certain rotation number
at infinity, where the monodromy is lifted to the universal cover of a
fibre (\ref{spinissquare}). 

\vspace{0.2cm}

\noindent Sections correspond to invariant zero-dimensional homology classes
of a certain type.  Here ``invariant'' means ``invariant with respect to the
monodromy diffeomorphisms which encode the topology of the fibration''.  
One can also ask about invariant one-dimensional homology 
or homotopy classes.  The monodromy
homomorphism

$$\rho: \pi_1 (\sss^2 \backslash \{ \Crit \}) \rightarrow \Gamma_g$$

\noindent has image the mapping class group $\Gamma_g$ of the generic smooth
fibre.  Composing with the natural map $\pi: \Gamma_g \rightarrow 
Sp_{2g} (\zz)$
from the mapping class group to the symplectic group, we obtain the 
cohomological monodromy which describes the action of the monodromy maps on
the symplectic vector space $(H^1 (\Fibre), \langle \cdot \rangle)$, where the
symplectic form is given by the cup-product.  It is an easy fact that there
is an identity $H^1 (X) = H^1 (F) ^{\rho}$ between the first cohomology of the
total space of the fibration and the invariant classes

$$H^1 (F)^{\rho} \ = \ \{ a \in H^1 (F, \zz) \ | \ f^* (a) = a \ \ \forall  f 
\in \mathrm{Image}(\rho) \}.$$

\noindent In particular, if $b_1 (X) > 0$ then $\pi \circ \rho$ has a trivial
subrepresentation.  Related properties of this cohomological monodromy
encode some of the fundamental properties of K\"ahler as opposed to 
symplectic four-manifolds: the Hard Lefschetz theorem states that if $X$ is
K\"ahler then this invariant subspace is a symplectic subspace, in particular
of even dimension.

\vspace{0.2cm}

\noindent The invariant cohomology is however at first glance quite 
insensitive - it can be defined from the data of the homology classes and
not the isotopy classes of the vanishing cycles.  One might expect to obtain
refinements of results such as the Hard Lefschetz theorem by passing from
invariant homology to invariant homotopy classes.  Our second result,
which first appeared in the author's D.Phil thesis \cite{ivanthesis}
and answers a
question due to Ludmil Katzarkov, 
is in that sense a disappointment:

\begin{Thm} \label{noinvariant}
Let $f: X \rightarrow \sss^2$ be any symplectic Lefschetz fibration;
write $F$ for 
a fixed smooth fibre and $\rho$ for the monodromy representation
as above.  Then there
is no finite set of disjointly embedded simple closed curves 
$\{ C_1, \ldots, C_r \} \subset F$ 
whose union is preserved to isotopy by $\rho$.  In particular, no free isotopy
class of embedded curve is invariant under the monodromy.
\end{Thm}

\noindent By way of consolation, we find that the vanishing cycles of
any Lefschetz fibration must fill the fibre $F$ (their complement is a bunch
of discs)\footnote{Apparently this result was classically known
  to Thurston and others, but I have been unable to trace a written
  account.}.   If the vanishing cycles did not fill, then by inserting 
a large genus surface with two boundary components in place of an annular
neighbourhood of a curve disjoint from the vanishing cycles, one could 
build a Lefschetz fibration with arbitrarily large fibre genus and fixed
number of critical fibres.  Thus one can deduce the ``one curve'' case of
(\ref{noinvariant}) from any genus-dependent bound on the minimal
length of positive
relations in mapping class groups;  however, such bounds have been
established only using gauge theory. 

\vspace{0.2cm}

\noindent Two well-known finiteness theorems in the geometry of
holomorphic fibrations are the Arakelov-Parsin theorem and Manin's
theorem (for a treatment from our perspective see for instance
\cite{Jost-Yau}).  The first has the
particular consequence that there are only 
finitely many holomorphic Lefschetz fibrations
with fixed fibre genus and fixed set of critical values; the second says
that any non-trivial holomorphic fibration admits at most finitely many
holomorphic sections.  The first of these results is false in the
symplectic setting;  there are infinitely many homeomorphism types of
symplectic manifold which admit Lefschetz fibrations by fixed genus
curves with a fixed number of critical values  (which can be placed at
any given set of points in the sphere).  Such examples are given in
\cite{ivanthesis} and this phenomenon is discussed in more detail in
\cite{Finsternparsin}.  An analogue of Manin's theorem, however, can be
deduced from (\ref{withthedisc}).

\begin{Cor}
Let $f: X \rightarrow \sss^2$ be a Lefschetz fibration.  Then only
finitely many  homotopy classes in $X$ contain  smooth sections of
$f$.
\end{Cor}

\noindent Note that this is false even for holomorphic fibre bundles
over bases of 
higher genus:  for instance, the trivial fibration $\rot^4 \rightarrow
\rot^2$ admits infinitely many distinct homotopy classes of smooth
section.  (It is also false in general at the level of homology; there
may be infinitely many homology classes with algebraic intersection
number one with the fibre.)  Similarly, one can contrast with the
infinite families of 
symplectic representatives for certain reducible homology classes -
multiples of a fibre of a Lefschetz fibration - 
given in \cite{ivansympsub}.  Any section of a Lefschetz fibration can
be made
a symplectic submanifold if we perturb the symplectic form by a
deformation equivalence, and our proof shall imply that the symplectic
isotopy class of a 
section is determined by its homotopy class.

\vspace{0.2cm}

\noindent \textbf{Acknowledgements:}  
Thanks to Ludmil Katzarkov and Gang Tian for motivating questions and
to Steve Kherchoff and Andras Stipsicz 
for valuable conversations.  Paul Seidel kindly identified some errors 
in an earlier draft of the paper. 

\section{Shearing the hyperbolic disc}

Our approach to (\ref{withthedisc}) will be to analyse the action of
the monodromy determining a Lefschetz fibration on the universal cover
of a fixed smooth fibre\footnote{Recall that we always assume that the
  fibre genus is at least
two.}.  From elementary covering space theory
we know that any homeomorphism of a topological space lifts to a homeomorphism
of the universal cover.  The lift is unique if we make the covering
map a map of pointed spaces by introducing base-points. For a
homeomorphism $\phi$ of a space $X$ we shall write $\tilde{\phi}$ for
any choice of lift to the universal cover.  The
uniformisation theorem allows us to view the two-dimensional fibre as
a quotient of the hyperbolic disc $D$ with its Poincar\'e metric by a
discrete subgroup of the isometry group of $D$.  Crucial to us will be
the circle at infinity $\circinf$ which naturally compactifies $D$.
Most of the material of the following section is standard and can be
found in work of Nielsen \cite{Nielsen};
see also the accounts of Thurston's work in \cite{Asterix} and
\cite{Thurstondynamics}.  For completeness we review the necessary
ideas.  Fix a 
covering projection $D \stackrel{\pi}{\longrightarrow} \Sigma_g$ and
an identification of $\pi_1 (\Sigma_g) \cong G$ with a discrete
subgroup $G < \mathrm{Isom}(D)$.  Write $e$ for the identity element
of $G$.

\vspace{0.2cm}

\noindent The Poincar\'e metric on the unit disc is entirely specified
by declaring that its isometry group consists of the M\"obius
transformations of the complex plane which preserve $D$.  With respect
to this metric, geodesics in $D$ are given by arcs of circles
which are orthogonal to the boundary $\circinf$.   A pair of distinct
points on $\circinf$ determines a unique
geodesic in $D$ and hence
geodesic $\gamma: (-\infty,\infty) \rightarrow \Sigma_g$ which may or
may not close.  Any choice of lift of a simple closed curve in
$\Sigma_g$ gives a curve
in $D$ 
meeting the boundary at two points which is
bounded isotopic to the geodesic corresponding to this boundary
pair. Since the group $G$ acts on $D$ preserving geodesics, it follows that
the action extends to an action on $\circinf$  by
orientation-preserving homeomorphisms.  For any $u \in G \backslash {e}$, the
induced map 
$u_*: \circinf \rightarrow \circinf$ has two fixed points, one
attracting and one repelling;  under iterations, all points on the
circle save for the repelling fixed point
are compressed near the attracting fixed point.  (Such a homeomorphism
of the disc and/or circle is called a \emph{hyperbolic} map.)  The geodesic
determined by the two fixed points projects to a closed curve
in $\Sigma_g$ which is the unique geodesic in the free homotopy class
associated to $u \in G$.  As we vary over $u \in G$ the fixed points
vary over a dense set of $\circinf$.  Note importantly that the
hyperbolic covering transformation has \emph{no fixed points} in the
interior of the disc $D$.  The following is our principal tool.

\begin{Prop} \label{wecanlift}
\begin{itemize}
\item Let $\tau: \Sigma_g \rightarrow \Sigma_g$ be any homeomorphism of the
surface.  Then each lift of $\tau$ to $D$ induces an
orientation-preserving homeomorphism of $\circinf$, and the set of
homeomorphisms one obtains in this way depends only
on the isotopy class of $\tau$ in the mapping class group.

\item The lifts of a Dehn twist about $C \subset \Sigma_g$ to $D$
  fall into two
  families:  those which fix some point inside the interior of $D$,
  which then fix a component of the complement in $D$ of a
  neighbourhood of the
  locus $\Lambda_C$ of all lifts of $C$, and those which fix no point of $\Int
  (D)$.  The members of the first family are permuted by the
  conjugation action of $G$ on the set of lifts.

\item Let $\delta: \Sigma_g \rightarrow \Sigma_g$ be a positive Dehn
  twist about some curve $C \subset \Sigma_g$.  Then any lift
  $\tilde{\delta}$ of $\delta$ to $D$ which has a fixed point in the
  interior of $D$ fixes a countable number of
  points of $\circinf$ and moves all other points in the same sense
  (clockwise in the convention in which holomorphic Lefschetz
  fibrations have positive Dehn twist monodromies).  
\end{itemize}
\end{Prop}

\begin{Pf}
The homeomorphism induces an automorphism of $G$
via $g \mapsto \tau g \tau^{-1}$ (which still covers the identity).
Using this we 
can induce an automorphism of the end-points of the geodesics
corresponding to elements of $G$;  since these form a dense subset of
the circle we therefore have a homeomorphism of $\circinf$.  It is
easy to check that this preserves the
cyclic ordering of the circle and is orientation-preserving. If we
change $\tau$ by isotopy we will not change the isotopy class of the
image of any curve under $\tau$; since there is a unique geodesic in
each isotopy class of simple closed curve, and it is the action on
these geodesics which determines the action on $\circinf$, the
independence of the extension under isotopies follows.

\vspace{0.2cm}

\noindent We give an explicit description of the lifts
which fix a region in the proof of the third part of the lemma;  it is
clear that these are permuted by the conjugation action of $G$.  Now suppose a
lift $\tilde{\delta}$ fixes some point of $\Int (D)$.  This lies
inside some region of 
$D \backslash \Lambda_C$ - it cannot lie inside $\Lambda_C$ since the
Dehn twist itself does not preserve any point of the curve $C$ (on
which it acts by the antipodal map).  For this particular
region  there is a unique lift $\tilde{\delta}_0$ of the twist which
fixes the open set given by the intersection of the region and the
exterior of a neighbourhood of the boundary lifts of $C$.  The product
$(\tilde{\delta})^{-1} \tilde{\delta}_0$ now covers the identity on $\Sigma_g$ and is
hence a deck transformation; but it also fixes a point of the
interior of $D$, by construction, which means that it must be the
trivial deck transformation $e \in G$.  It follows that a lift of the
Dehn twist fixing some interior point necessarily co-incides with the
lift fixing the region containing the point.

\vspace{0.2cm}

\noindent The
homeomorphism $\tau_C$ given by twisting about $C$ is supported inside
an annular neighbourhood $A$ of $C$;  standard covering space theory
implies that it can be lifted to be the identity
in a  single component $K$ of the universal
cover of $\Sigma \backslash A$.  Fix a base-point $P$
inside the region $K$.
Given the definition of the action of a mapping
class on the 
  circle at infinity, we may understand the effect of the Dehn twist
  about a curve $C$ as follows.  Choose a point $p$ on
  $\circinf$; visualise $p$ as lying at twelve o'clock.  This point,
  together with the base-point $P$,
  determines a quasigeodesic $L: [0,\infty) \rightarrow D$ and a quasigeodesic
  $l$ in
  the surface $\Sigma_g$ which may intersect $C$, perhaps infinitely
  often.  The effect of the Dehn twist $l \mapsto \tau_C (l)$ is to
  insert a copy of $C$ at
  each of these intersection points.  Start at the image of $P$ in
  $\Sigma_g$ and
  consider moving along the new curve $\tau_C (l)$; then upstairs
  we start at $P$ and move along $L$ until we meet some lift of $C$,
  then we turn right and move along
  this lift of $C$ until we meet the adjacent lift of $l$, then we
  turn left and move in the initial direction along this new parallel
  lift of $l$ 
  until we meet the next intersection point with a lift of $C$, turn
  right and so 
  forth (draw a picture!).  This process eventually
  meanders into $\circinf$, by
  definition giving rise to the image of  $p$ under the
  homeomorphism of $\circinf$ induced by this particular lift of the
  Dehn twist about $C$. Suppose for contradiction this image point
  $\tilde{\tau}_C (p)$ lies 
  to the left of 
  $p$ (as viewed from inside the disc, i.e. anticlockwise along the circle
  from $p$).  Then after some finite number of intersections with $C$
  and zig-zags in the above procedure, the extension of the lift of
  $l$ on which we are travelling also has end-point anticlockwise from
  $p$.  Suppose for simplicity that this happens in fact after just
  one zig-zag, that is that $C$ and $l$ meet only once in $\Sigma_g$.
  Then we have a geodesic triangle inside $D$ comprising the arc of
  $C$ and the two arcs of lifts of $l$, which by our assumptions must
  meet at some third point. But now two of the angles of our triangle (at the
  intersection points of the arc of $C$ and of lifts of $l$) sum
  to $\pi$; this is impossible in hyperbolic trigonometry.  

\vspace{0.2cm}

\noindent For a
  larger number of zig-zags one can complete by induction, the crucial
  point being that the intersection angles of lifts of $C$ and of $l$
  are always equal in the disc.  It follows that all points that are
  moved at all are indeed moved in one particular sense, determined by
  the positivity of the Dehn twist homeomorphism downstairs.  Moreover
  it follows that some points of the circle are moved by the lift
  $\tilde{\tau}_C$, coming from closed geodesics which do
  intersect $C$.

\vspace{0.2cm}

\noindent Finally, to see that countably many points are fixed,
consider again the fixed lift $K$
in $D$ of the 
complement $\Sigma \backslash A$.  The boundary of this
complement downstairs lifts to give 
pairs of curves with equal endpoints which are endpoints of lifts of
the geodesic in the free homotopy class of $C$.  There are countably many
of these in the boundary of $K$ (of course the endpoints of all the
lifts of $C$ give a dense subset of $\circinf$ but this is not true if
we restrict to those at the boundary of the one region).  Now the
homeomorphism $\tau_C$ given by twisting about $C$ is the identity
inside $K$ and hence preserves (setwise!) each of the boundary curves
of $K$.  It follows that the induced map on $\circinf$ has a countable set
of fixed points.  
\end{Pf}

\begin{Rmk}
Note that it is \emph{not} true that every lift of a Dehn twist has
this nice behaviour at the circle at infinity.  For if we compose a
lift with an interior fixed point with some large power of a fixed
hyperbolic deck transformation, the compression of $\circinf$ near the
attracting fixed point of the hyperbolic will destroy the clockwise
shear.  (This example was pointed out to me by Paul Seidel.)
Fortuitously the applications we have in mind necessitate choosing
lifts of Dehn twists which \emph{do} have interior fixed points.
\end{Rmk}

\noindent So far the non-uniqueness of lifts presents a problem for
applications.  If we suppose that our fibration has a section,
however, this problem is entirely removed.  Recall that if we fix a
base-point $\tilde{p}$ in $D$ covering a base-point $p$ in $\Sigma_g$
any homeomorphism $\tau$ of the base has a distinguished lift
$\tilde{\tau}_{(p)}$ to $D$.  If the
homeomorphism fixes $p$ then the lift fixes the distinguished point
$\tilde{p}$ (but in general induces a non-trivial permutation of the
other lifts  of $p$).  A section of a Lefschetz fibration precisely
amounts to a point $p \in \Sigma_g$ fixed by each of the Dehn twist
monodromies $\delta_i$.  In fact more is true; the point is then
fixed by the isotopy from the product homeomorphism $\prod \delta_i$
to the identity. If we lift each of the Dehn twists in turn, we
induce a sequence of homeomorphisms of 
the circle whose product must be the identity.  For the final isotopy
does not have any effect on $\circinf$ by the last part of
(\ref{wecanlift}), and so the product of the lifts $\prod
\tilde{\delta_i}$ acting on $\circinf$ must be the same as the lift of
the identity, which is the result of that isotopy downstairs in
$\Sigma_g$.  But the whole process has fixed the point $\tilde{p}$ and
only the identity deck transformation fixes a point of $D$.  Thus the
final homeomorphism on $\circinf$ co-incides with the extension to
$\circinf$ of the lift of the identity of $\Sigma_g$, which is clearly
trivial.  However, by lifting the various Dehn twist homeomorphisms in
turn we have collected another piece of information;  a winding
number, telling us how often we have spun the circle.  (This winding
number can be viewed as a translation number if we lift now to the
universal cover $\rr$ of the circle $\circinf$ but the meaning should
be clear without this additional layer of notation.)

\begin{Lem} \label{spinissquare}
The sequence of homeomorphisms $\tilde{\delta_i}$ 
rotates the circle clockwise by $2\pi k$, where the positive
integer $k = -s\cdot s$ is given by the negative of the
self-intersection of the section $s$  of
the Lefschetz fibration defined by the point $p$.
\end{Lem}

\begin{Pf}
Fix a small disc centred on $\tilde{p}$. Since $\tilde{p}$ is fixed by
all the homeomorphisms of the disc $\tilde{\delta_i}$ and the final
isotopy, this disc is mapped eventually to another topological disc
enclosing $\tilde{p}$.  It follows that the boundary of a small disc
around $\tilde{p}$ is to homotopy rotated a certain number of times by
the product of the homeomorphisms;  if we now take two points
$\tilde{p}, \tilde{q}$ which are nearby lifts of close points $p,q$ in
$\Sigma_g$ then under the sequence of homeomorphisms the point
$\tilde{q}$ will move around $\tilde{p}$ a number $m$ times.  This
number is precisely the number of intersection points of the sections
defined by $p,q$ respectively, which is just the self-intersection
$|s\cdot s|$.  On the other hand, by considering a radial projection
along geodesic arcs,
since the homeomorphisms of the disc preserve the cyclic ordering on
$\circinf$ we see that this number is the same as the rotation number
at the boundary of the disc.  In more formal language, there is a
short exact sequence

$$0 \rightarrow \zz \rightarrow \Gamma_{g,1} \rightarrow \Gamma_g ^1
\rightarrow 0$$

\noindent where the middle term is the mapping class group for a
surface with one boundary circle and the right hand term the mapping
class group for a punctured surface.  The copy of $\zz$ is generated
by the Dehn twist about the puncture;  our Dehn twist monodromies
naturally lift from $\Gamma_g ^1$ to $\Gamma_{g,1}$ but their product
may differ from the identity in $\Gamma_{g,1}$ by some integer in the
kernel of the natural map $\Gamma_{g,1} \rightarrow \Gamma_g ^1$, and
this is the self-intersection of our section.  There is another short
exact sequence

$$1 \rightarrow \zz \rightarrow \Homeo^{\flat}(\rr) \rightarrow
\Homeo(\circinf) \rightarrow 1$$

\noindent where all homeomorphisms are orientation-preserving and the
middle group denotes the homeomorphisms of the real line $f$ for which
$f(x) + 2\pi = f(x+2\pi)$.  The copy of $\zz$ is generated by
translations.  Then the statement of the lemma amounts to saying that
the first short exact sequence is induced from the second under the
representation $\Gamma_g ^1 \rightarrow \Homeo(\sss^1)$.

\vspace{0.2cm}

\noindent Being explicit again, note
from the description with $p,q$ that the sign of the self-intersection
of any section is determined by the sense of the
rotations on $\circinf$.  It follows that the sign of the
self-intersection of the 
section is always the same for any Lefschetz fibration determined by
positive Dehn twists;  but we know of many fibrations containing
sections of negative square, from blowing up base-points.
\end{Pf}

\section{Irreducible Lefschetz fibrations}

Since one can always introduce critical fibres in Lefschetz fibrations by
blowing up points, it is standard to assume that all Lefschetz fibrations have
no spherical components in fibres.  The genus zero fibrations are then 
Hirzebruch surfaces whilst the genus one fibrations are elliptic surfaces
\cite{BPV}.  All the fibrations are holomorphic, and moreover all are given
by fibre summing certain ``irreducible'' building blocks by identity
diffeomorphisms.  This prompted the following question:

\begin{Question}[Gang Tian]  Is every simply-connected Lefschetz fibration
a fibre sum of holomorphic Lefschetz fibrations?
\end{Question}

\noindent From the known restrictions on the minimal length of mapping
class group words, one can find examples showing that the assumption of
simple connectivity is necessary \cite{ivanthesis}.  Our next example
shows that it is
not sufficient:

\begin{Prop} Any pencil of curves of genus $g > 5$ on a symplectic
  non-K\"ahler $K3$ surface gives rise to a Lefschetz fibration which
  is not a K\"ahler sum.
\end{Prop}

\noindent There are many - gauge-theory based - constructions of such
  fake $K3$ surfaces
  \cite{FinStern};  the Lefschetz pencils, and hence fibrations on
  blowing up base-points, are provided by Donaldson's existence
  result.  The fibres of the pencil may have arbitrarily large genus
  by scaling the symplectic form and computing with the adjunction formula.
\vspace{0.2cm}

\begin{Pf}
Let $Z' \rightarrow \pp^1$ be the fibration constructed from our fake
$K3$ surface.  Note that $\sigma + e$ - the sum of signature and
Euler characteristic - is invariant under blowing up and down.  If we
write, for contradiction, that $Z' = W_1 \sharp_{\Fibre} W_2$ is a
fibre sum (of non-trivial Lefschetz fibrations, that is fibrations
with non-zero number of singular fibres), then by some easy
computations we have 

\begin{Eqn} \label{noughtnote}
\sigma(Z') = \sigma(W_1) + \sigma(W_2); \ e(Z') = e(W_1) +
e(W_2) - 2e(F)
\end{Eqn}

\noindent where $F$ denotes a smooth fibre.  It follows, in obvious
notation, that

\begin{Eqn} \label{firstnote}
(\sigma_1 + e_1) + (\sigma_2 + e_2) = 8+2e(F) \ < \ 0 
\end{Eqn}

\noindent and hence for $W_1$ say we have $\sigma_1 + e_1 < 0$.
But this forces, by the classification of complex surfaces, the
manifold $W_1$ to be the blow-up of an irrational ruled manifold
$\Sigma_h \tilde{\times} \sss^2$ for some $h$.  Since $\pi_1 (W_1)$ is
$\pi_1 (\Sigma_h)$ and also a quotient of $\pi_1 (\Sigma_{g(F)})$ by
vanishing cycles we see $h \leq g(F)$.  Indeed if $h=g(F)$ then all
the vanishing cycles for $W_1$ are nullhomotopic and it is a trivial
piece in the decomposition, a contradiction.  Hence $h < g(F)$.

\vspace{0.2cm}

\noindent We claim that in fact $2h \leq g(F)$.  To see this, take a
basis for the homology of the irrational ruled manifold comprising a
section $[s]$ and fibre $[f]$ (the latter is canonically defined, the
former only up to twisting $[s] \mapsto [s\pm nf]$).  By the
adjunction formula, note that any complex curve representing any
section of the fibration $\Sigma_h \tilde{\times} \sss^2 \rightarrow
\Sigma_h$ has genus $g_{[s\pm nf]} = h$.  Since $g(F) >h$ we must have
$F = a[s] + b[f]$ with $a \neq 1$, and hence $a \geq 2$ since the
symplectic form is positive on $F$.  Finally, McDuff \cite{McD-S} has
shown that
for a symplectic form $\omega= P.D(a[s]+b[f])$ on an irrational ruled
manifold $W$ with fibre $\sss$ we have

$$\omega^2 (W) > \omega(\sss)^2 \ \Rightarrow a^2 [s]^2 + 2ab > a^2.$$

\noindent Applying adjunction to $F = a[s] + b[f]$ shows that 

$$2g(F) - 2 = aK_W \cdot [s] - 2b + a^2 [s]^2 + 2ab = a(2h-2) +
(a-1)(a[s]^2 + 2b);$$

\noindent combining this with McDuff's result and $a \geq 2$ gives
$g(F) \geq 2h$ as claimed.

\vspace{0.2cm}

\noindent Now $\sigma_1 + e_1$ (invariant under blow-ups) is equal
to $2 e(\Sigma_h)$ whilst $2e(\Sigma_h) - 2e(F) > 8$ from the
original stipulation that $g(F) > 5$.  It therefore follows from
(\ref{firstnote}) that $\sigma_2 + e_2 < 0$ and hence that $W_2$ is
also a blow-up of an irrational ruled manifold.  From here the readers
should be able to find their own conclusion to the proof; the combined
restriction of having both pieces of a decomposition (blow-ups of)
irrational ruled manifolds is too considerable.
\end{Pf}

\noindent The details of the proof are surprisingly technical.  In
fact there is a much more elementary obstruction at work.
The Lefschetz fibrations produced by Donaldson's construction are
\emph{never} fibre summations.  This was first proven
by Stipsicz \cite{Stipsicz} using Seiberg-Witten theory, and in
particular the fact that  a symplectic manifold with $b_+ > 1$ cannot
contain a symplectic sphere of
positive square.  The above example (historically the first) can be seen as a 
particular realisation
of those kinds of method.  In fact the general case follows easily from our 
remarks on the hyperbolic
disc above.  The proof
will recover, as an addendum, the fact that symplectic spheres which arise as
\emph{sections} of Lefschetz fibrations
can never have positive square.

\begin{Prop}
\begin{itemize}
\item No Lefschetz fibration can
contain a section of positive square.
\item If a Lefschetz fibration contains a section of square zero, it is a
trivial product $\Sigma_g \times \sss^2$.  
\item If there is a section of square
$-1$ then the fibration cannot split as a non-trivial fibre sum.
\end{itemize}
\end{Prop}

\begin{Pf}
The first two statements are immediate from the last remarks in the proof of
(\ref{spinissquare}).  Suppose then for
contradiction that the Lefschetz fibration $X \rightarrow
\pp^1$ contains a section of square $(-1)$ and splits as a non-trivial
fibre sum $X = W \sharp_{\Fibre} W'$.  Choose a positive relation

$$\delta_1 \ldots \delta_r \, \gamma_1 \ldots \gamma_s = 1$$ 

\noindent for $X$ which
is a factorised product of relations for the constituent pieces $W,
W'$.  Fix a covering projection $\pi: D \rightarrow \Fibre$ as above
and a point $\tilde{p}$ in $D$ which projects to the intersection of
the section with the distinguished fibre.  We have a canonical
sequence of lifts of monodromies $\tilde{\delta_i}, \tilde{\gamma_j}$
whose product induces the single full rotation of $\circinf$ by
$2\pi$.  On the other hand, once we have lifted the first $r$
monodromies, we have lifted a word which is the identity, and hence we
must induce a hyperbolic covering transformation on $\circinf$ by the
arguments as above.  Thus the sequence of lifts serves to factorise
the $2 \pi$-rotation of $\circinf$ as a product of two hyperbolic
automorphisms, each of which is a factorised product of clockwise shears. 

\vspace{0.2cm}

\noindent But this is impossible, from our description of these
hyperbolic elements.  The first hyperbolic element has two fixed
points.  Therefore under the factorisation into twists, these points
have either made at least one rotation of the circle or they have been
fixed throughout.  In the former case, points on one side of the
repelling fixed point must have moved by more than one full rotation
of the circle, since the cyclic ordering is preserved after each
twist.  But then since the total translation is increasing, composing
with the twists arising from $W'$ cannot yield a total translation of
one $2\pi$ rotation.  So suppose that the two fixed points have been
fixed throughout all the Dehn twist lifts for $W$.  They must then be
moved by one rotation by the lifts for $W'$.  However, by the same
argument, the two fixed points for $W'$ are now fixed by all these
latter rotations, and so the fixed points for $W$ cannot pass them in order
to move around the boundary of the disc.
\end{Pf}

\noindent We should make one point clear: it is not the case that all
sections of a given fibration must have the same square.
K\"ahler genus two fibrations with no reducible fibres were classified
by Chakiris (cf. \cite{ivanhodge}).  From that
classification, it follows that the fibre sum (always by the identity
diffeomorphism) of four copies of the fibration associated to the
genus two pencil on $K3$ is isomorphic to the fibre sum of six copies
of the fibration on the rational manifold $(\sss^2 \times \sss^2) \sharp 12
\overline{\pp}^2$; they each have full monodromy group and $120$
singular fibres.  On the other hand, in these decompositions the
holomorphic fibration is shown off with sections of square $-4$ and $-6$
respectively.  There is no analogue of this phenomenon for
elliptic fibrations, where the adjunction formula fixes the square of
any section.

\section{Invariant homology is not homotopic}

The most classical way to detect that a symplectic manifold is not
K\"ahler is to
show that the first Betti number is odd.  Recall that the first
homology of a Lefschetz fibration comprises the homology classes in
the fibre which are invariant under the monodromy action:  the
representation of $\pi_1 (\sss^2 \backslash \{\Crit \} )$ on $H_1
(\Fibre)$ has a trivial subrepresentation $H_1 (X)$.  It has become
clear that treating the vanishing cycles up to isotopy and not just
homology gives access to additional information on the underlying
manifold.  This was one motivation for the following

\begin{Question}[Ludmil Katzarkov] 
Can there be a simple closed curve in a fibre of a Lefschetz fibration
which is invariant (up to isotopy) under all the monodromy
diffeomorphisms of the fibration?
\end{Question}

\noindent Katzarkov \emph{et al} \cite{KPS} gave a negative
answer, based on Hodge theory, 
for K\"ahler fibrations arising from pencils of high degree on
surfaces with vanishing first Betti number.  The methods of
the previous  section make short
work of the more general case.

\begin{Prop}
There is no finite set of simple closed curves with embedded union
which is preserved to isotopy by the monodromy of any Lefschetz
fibration.
\end{Prop}

\begin{Pf}
We treat the case of a single curve first.  Suppose $C \subset
\Sigma_g$ is such a curve.  Fix a lift of $C$ to the disc extended to
the boundary $\circinf$ giving two distinguished points $a, b \in
\circinf$.  Since all the monodromy diffeomorphisms fix $C$ and in
particular all the $\delta_i$ fix $C$, we can choose lifts of each
$\delta_i$ in turn which fix \emph{the whole lift} $\tilde{C}$.  To
see this, note that we can always lift a homeomorphism to be the
identity in one lift of the complement of its support, and here we
choose all these lifts to contain $\tilde{C}$.  It follows that under
the various monodromies upstairs, $a,b$ are always fixed.  However,
the product of these monodromies must be a deck transformation on
$\circinf$.  This is impossible;  the two fixed points $a,b$ must be the
fixed points of the hyperbolic, but because all the Dehn twists
have moved points in a single sense, we see that in small arcs
either side of say $a$ the points of $\circinf$ are all moved
clockwise.  This is impossible if $a$ is either an attracting or
repelling fixed point.

\vspace{0.2cm}

\noindent To see that we cannot permute a collection of simple closed
curves $\gamma_1, \ldots, \gamma_q$ we use induction, with the above
being the base step.   Fix a
distinguished curve $\gamma_1$; by the above there must be some vanishing
cycle $\psi$ which has non-trivial geometric intersection number with
$\gamma_1$.  Consider the curves $\gamma_1 ^{(n)}$
given by applying the $n$-th power of the Dehn twist about $\psi$ to
$\gamma_1$.  Eventually we must have some identity $\gamma_1 ^{(n)} =
\gamma_1$ for some $n$.  But then the group generated by the Dehn
twists about $\gamma_1, \psi$ has a relation:  the twist about
$\gamma_1$ and the $n$-th power of the twist about $\psi$ commute.  It
is known however that the group generated by the Dehn
twists about two curves with non-trivial geometric intersection number
is free of rank two if the intersection number is greater than one
and the braid group on three strings if the number is precisely one.
No relation of the given form can exist in either of these cases.
\end{Pf}

\noindent We can put a different interpretation on this last result.
For a simple closed curve to be invariant under the Dehn twists about
each vanishing cycle it would have to be disjoint from all of these;
conversely, anything disjoint would certainly be invariant. 

\begin{Cor}
For any Lefschetz fibration, the vanishing cycles fill up the fibre;
their complement is a bunch of discs.
\end{Cor}

\vspace{0.2cm}

\noindent It follows that the pattern of vanishing cycles on the fibre
$\Sigma_g$ amounts (in a generic and transverse picture) to a certain
four-valent graph on the surface which
is the one-skeleton of a cell decomposition.  One can show that one
half the
number of faces in this decomposition gives an upper bound on the
number of distinct homotopy classes of section of the fibration, but
in general the bounds one gets this way are not remotely sharp.  More
curious is the following observation. So far the particular hyperbolic
structures on surfaces that we have considered have been irrelevant
and certainly not canonical.  However, Kherchoff showed in
\cite{Kherchoff} that the geodesic length function is convex along
earthquake paths.  It followed that for any collection of simple
closed curves which fills a surface, there is a distinguished element
of Teichm\"uller space $T_g$ on which the total length of the
geodesic representatives of these curves is minimised.  Counting
repeated isotopy classes of vanishing cycle with multiplicity, this
total length defines a value $l(V)$ for any collection of vanishing
cycles defining 
a Lefschetz fibration.  Such collections are indexed by a braid group
$\Br_n$ as explained in \cite{ivanthesis}.

\begin{Defn}
Let $f:X \rightarrow \sss^2$ be a genus $g$ Lefschetz fibration defined by a
family of positive relations $I(\Br_n)$ indexed by a suitable braid group.
The \emph{length} $l(X,f)$ is
the minimum $\min \{ l(V_{\iota}): \iota \in I(\Br_n) \}$ and is an
invariant of 
the Lefschetz fibration.
\end{Defn}

\noindent Such invariants are
unlikely to have subtle geometric content, since they do not
depend on the \emph{order} of the Dehn twist monodromies, and only on
the supporting curves.  Nonetheless it would be interesting to
understand precisely what if anything they do capture.  For instance,
the length may plausibly be related to the minimal possible value of the total
number of geometric intersections of the vanishing cycles, as we vary
over the braid group.  This latter invariant is related to  Floer theory,
and the sum of the
ranks of the Floer homologies over all pairs of distinct vanishing cycles in
the fibre.


\section{Finitely many sections}

\noindent Manin's theorem \cite{Jost-Yau} asserts that a
non-isotrivial holomorphic
fibration has only finitely many holomorphic sections.  We have the
following analogue for Lefschetz fibrations:

\begin{Thm}
Let $f: X \rightarrow \sss^2$ be a symplectic Lefschetz fibration.
Then $f$ admits only finitely many homotopy classes of geometric
section.
\end{Thm}

\begin{Pf}
Any section can be trivialised over a large disc in the base, containing
the critical values, and viewed as a point $p$ in a fixed fibre
$\Sigma$ which is 
invariant under all the monodromy maps.  That is, $p$ is disjoint from the
support of all the Dehn twist monodromies.  The product of these is
isotopic to the identity; under the isotopy, which determines how to close
the fibration over the disc to one over the sphere, $p$ moves through some
loop $u \in \pi_1 (\Sigma)$.  There is a section through $p$ iff this loop is
nullhomotopic in $\Sigma$.

\vspace{0.2cm}
 
\noindent It follows that if we perturb $p$ inside a component $D$ of $\Sigma
\backslash \{ C_i \}$, where the $C_i$ are the vanishing
cycles of the fibration transported to lie in the fixed fibre,  then
the loop $u$ associated to $p$ is unchanged.  From this it is easy to
see that if there is a section through $p$ its homotopy class is
independent of the choice of $p$ in $D$.  So the number of homotopy classes of
section is bounded above by the number of cells in the cell decomposition
of $\Sigma$ defined by the vanishing cycles, and is in particular
finite.
\end{Pf}

\noindent One can improve the actual bound on the number of homotopy
classes of sections, but the results seem a long way from being sharp
in examples.  If one is given a component $D$ for which there is a
section of the fibration, it is in principle mechanical to determine
which other regions of $\Sigma \backslash \{ C_i \}$ admit sections.
For taking a path between points in the two regions, and the image of
the path under all the monodromy twists, we have a loop in $\Sigma$
which should be trivial in $\pi_1 (\Sigma)$.  Such a relation amongst
the vanishing cycles is equivalent to the existence of a suitable
closed polygon in the hyperbolic disc.  Note
finally that the \emph{existence} question for sections of
Lefschetz fibrations remains open.


\bibliographystyle{amsplain}
\bibliography{main}

\end{document}